\documentclass[12pt]{article}
\usepackage{amssymb}

\usepackage{amsfonts,amssymb,eucal,amsmath}
\usepackage{graphicx}
\title{A New Fixed Point Theorem for Non-expansive Mappings and Its Application}
\author{{\small\sc Chunyan Yang}\\
\small {\it Department of Mathematics, Sichuan University,
         Chengdu 610064,P.R.China}\\
       \small  e-mail: 787090614@qq.com}
\date{\small Feb,12,2012}

\begin{document}
\maketitle

\textbf{Abstract.} { We use $KKM$ theorem to prove the existence
of a new fixed point theorem for non-expansive mapping:Let M be a
bounded closed convex subset of Hilbert space H, and
$A:M\rightarrow M$ be a non-expansive mapping, then exists a fixed
point of A in M, we also apply this Theorem to study the solution
for an integral equation,we can weak some conditions comparing
with Banach's contraction mapping principle.}



\section{\textbf{Introduction and Main Results}}\label{sect 1}

 \paragraph{    }It is well-known that $Banach$, $Schauder$,$Brouwer$  presented three kinds of fixed point theorems in 1910's-1930's. $Banach's$ fixed point theorem is also called contraction mapping principle , and  $Schauder's$ fixed point theorem is a generalization of $Brouwer's$ from finite dimension into infinite dimension. Specifically, for a compact operator A, if its domain is a  bounded non-empty closed convex subset, then there is at least one fixed point u : $A \cdot u = u$. All the above three  fixed point theorems need good operators, for example, the compactness for the operator can make the unit closed sphere in infinite dimensional space be a compact set.
 \paragraph{   } Here, a nature question is whether the fixed point theorem of other good operator or mapping could be established.
 \paragraph{  }A problem is whether there is a fixed point theorem for A when $k = 1$ (where A is called the non-expansive mapping). For solving this problem, the main tool we will use is $KKM$ theorem.It's well known that the Polish mathematician $Knaster$, $Kuratowski$, $Mazurkiewicz$ got $KKM$ theorem in 1929, they also applied it to the proof of $Brouwer$ fixed point theorem.
\paragraph{} \textbf{ Theorem1.1.}( Banach fixed point theorem(\cite{1})):
 Let $(X,d)$ be a complete metric space, and
  $A: X\rightarrow X$ be a contraction mapping, that is ,
  for any x,y $\in$ X, there holds $d(Ax, Ay)\leq k d(x, y), (0 \leq k<1)$.
  Then there is a unique fixed point $x*$ for $A$,  i.e. $Ax*=x*$.
 \\Banach's fixed point theorem can be applied to many fields in mathematics,especially to integral equation:
$$u(x)=\lambda\int_{a}^{b}F(x,y,u(y))dy+f(x),\ \ \ \ \ a\leq x \leq
b\eqno(1.1)$$
\paragraph{}  \textbf{ Theorem1.2.(\cite{2})} Assume

(a) $f :[a,b]\rightarrow R $ is continuous

(b) $F(x,y,u(y))=K(x,y)u(y)$,and $K :[a,b]\times[a,b]\rightarrow R $
is continuous.

Let $\Gamma=\max_{a\leq x,y\leq b}|K(x,y)|$,and

(c) real number $\lambda$ be given such that
$(b-a)\cdot|\lambda|\cdot\Gamma<1$.

Then the problem (1.1) has a unique solution $u\in X=C[a,b]$.\\
In this paper,we get:

\paragraph{}  \textbf{ Theorem1.3.} Let M be a bounded closed convex subset
of Hilbert space H, and $A:M\rightarrow M$ be a non-expansive
mapping, then exists a fixed point of A in M.
\paragraph{} Notice that our fixed point theorem is  in Hilbert Space , since any
linear continuous functional in Hilbert space can be represented
by the inner product, which is called $Riesz$ representation
theorem, while in the Banach space no such theorem.

We apply the above theorem to linear integral equation,we obtain
 \paragraph{}  \textbf{ Theorem1.4.} For  the integral equation(1.1)
 with $F(x,y,u(y))= K(x,y)u(y),$
 assume

(1)The function $f:[a,b]\rightarrow R$ is in $L^{2}[a,b].$

(2)$K(x,y)\in L^{2}([a,b]\times[a,b])$,$K(x,\cdot)\in
L^{2}[a,b]$,for any given $x\in[a,b]$

(3)Let the real number $\lambda$ be given such that,
 $|\lambda|\int_{a}^{b}\int_{a}^{b}|K(x,y)|^{2}dx dy \leq 1,if f\equiv0;$ and,
$$|\lambda|\int_{a}^{b}\int_{a}^{b}|K(x,y)|^{2}dx dy<1,if f\neq0 .$$
 \\Then (1.1) has at least one solution $u \in L^{2}[a,b].$

 \section{\textbf{Related Definitions and Well-known Theorems}}\label{sect 2}

 \paragraph{} \textbf{  Definition 2.1.} ($KKM$ mapping(\cite{3})): The mapping $G: X \rightarrow 2^{R^n}$ , where $X$ is a non-empty subset of the vector space
 $R^n$ , is called $KMM$ mapping. if $ \forall \{x_0, \ldots ,x_n\} \subset X$, we  have $Co(\{x_0, \ldots , x_n\})\subset \cup^n_{i=1} G(x_i) $,
 where $Co(\{x_0 ,\ldots ,x_n\})$ is the closed convex hull of $x_0 ,\ldots, x_n.$
 \paragraph{} \textbf{Theorem 2.1.} ($KKM$ Theorem (\cite{4}))
 Let $\Delta_n$ be n dimensional simplex which has vertices $\{e_0,\ldots ,e_n\}$
 and is in $R^{n+1}$.
  Let $M_0,\ldots ,M_n$ be n dimensional closed subsets in $R^{n+1}$,
  which satisfy $\forall \{e_{i0},\ldots ,e_{ik} \} \subset \{e_0,\ldots, e_n\}$,
  $Co(\{e_{i0},\ldots, e_{ik}\})\subset \cup_{j=0}^kM_{ij}$. Then $\cap M_i$ is not empty.

 \paragraph{}   \textbf{Theorem 2.2.} ($FKKM$ (\cite{5})): Let $F^n$ be Hausdorff linear topological space, X be a non-empty subset of $F^n$. Let $G: X\longrightarrow 2^{F^n}$ be a $KKM$ mapping, and $G(x)$ is weakly compact for any $x\in X$ . Then $\cap_{x\in X} G(x)$ is non-empty.
\paragraph{} \textbf{Definition 2.2.}( Semi-continuous(\cite{6})): Let $X$ be a Banach space, $X^*$ be its dual space. Let $T: X\rightarrow X^*$. If $\forall y \in X $, $\forall t_n \geq 0 $and $x_0+t_n\cdot y \in X$, we have $lim_{x_0+t_n\cdot y \rightarrow x_0}T(x_0+t_n\cdot y)\rightarrow T(x_0)$, we call T is $semi-continuous$ at point $x_0$. Furthermore if T is $semi-continuous$ at every point x in X, T is called $semi-continuous$ on X.
\paragraph{} \textbf{Definition 2.3.}( Monotonicity(\cite{7})): Let X be a Banach space, $X^*$ be its dual space. We call $T:\rightarrow 2^{X^*}$  monotonous , if $<u-v,x-y> \geq 0$ ,for $\forall x,y \in X$, $u \in T(x), v\in T(y)$.

\section{\textbf{Some Lemmas for the Proof of Our Theorem 1.3.}}\label{sect 3}
\paragraph{} We require the following Lemmas for proving the fixed point theorem of the non-expansive mapping using  the $KKM$ theorem .

\paragraph{} \textbf{Lemma 3.1}: Let M be a bounded closed convex subset of Hilbert space H, $L:M\rightarrow H$ be monotonous semi-continuous mapping. If $x_0\in M$, then
$$ <L(x_0),y-x_0> \geq 0, \forall y\in M ,$$
if and only if $$<L(y),y-x_0> \geq 0, \forall y\in M.$$
\\\textbf{Proof} Let $x_0 \in M$ such that :
$<L(x_0),y-x_0> \geq 0,\forall y \in M$. Since the monotonicity of
L, we have $$<L(y)-L(x_0),y-x_0>\geq 0\Leftrightarrow
<L(y),y-x_0>-<L(x_0),y-x_0>\geq 0,\forall y \in M.$$
 Thus ,
$$<L(y),y-x_0>\geq <L(x_0),y-x_0>\geq 0, \forall y\in M.$$
On the other hand,if $x_0 \in M$ then $$<L(y),y-x_0>\geq 0, \forall
y\in M.$$ We let $$y=h\cdot\mu +(1-h)x_0, \forall h \in
(0,1],\forall \mu \in M ,$$
since M is convex, then
 $y\in M$. Thus,
$$<L(y),y-x_0>=<L(x_0+h\cdot(\mu-x_0)),h\cdot(\mu-x_0)>\geq0(h>0).$$
We let $ h\rightarrow 0 $, then $x_0 + h\cdot(\mu-x_0)\rightarrow
x_0$, since L is semi-continuous , so $$<L(x_0),(\mu-x_0)> \geq
0,\forall \mu \in M.$$

\paragraph{} \textbf{Lemma 3.2}(\cite{8}): Let M be a bounded closed convex subset of Hilbert space H, and$L:M\rightarrow H, P:M\rightarrow 2^M$, and $P(y)=\{x\in M, <L(x),x-y>\leq0\},\forall y \in M$, then
 P is a $KKM$ mapping.
\\\textbf{Proof}: If P is not a $KKM$ mapping. By the definition of $KKM$ mapping,
$\exists \{x_0\ldots x_n\}=w$, $w\subset M$, let $\bar{x}=\sum
h_i\cdot x_i \in Co(w),(\sum h_i=1, h_i\geq0)$, then
 $\bar{x} \not\in \cup_{i=0}^n P(x_i)$, that is to say,
 $<L(\bar{x}),\bar{x}-y>>0$.
As the arbitrariness of $y$, let $y=\sum h_i\cdot x_i $, we
have$<L(\bar{x}),\bar{x}-x_i>>0$, then
$<L(\bar{x}),h_i\cdot(\bar{x}-x_i)>>0(\exists h_i>0)$, thus,
$$<L(\bar{x}),\sum h_i\cdot(\bar{x}-x_i)>=<L(\bar{x}),\sum h_i\cdot
\bar{x}-\sum h_i\cdot x_i>=<L(\bar{x}),\bar{x}-y>=0.$$
This is a
contradiction with $<L(\bar{x}),\bar{x}-y>>0,\forall y\in M$. Hence
P is a $KKM$ mapping.

\paragraph{}\textbf{ Lemma 3.3}: Let M be a bounded closed convex subset of Hilbert space H, $L:M\rightarrow H$ be a semi-continuous monotonous mapping, then $\exists x_0 \in M$ s.t. \[ <L(x_0),x_0-y>\leq0,\forall y\in M \]
\textbf{Proof}: Let $P:M\rightarrow 2^M$ and $P(y)=\{x\in M,
<L(x),x-y>\leq 0\},\forall y\in M$. According to Lemma 3.2,P is a
$KKM$ mapping. Mean-while, let mapping$J:M\rightarrow 2^M$ and
$$J(y)=\{x\in M,<L(y),x-y>\leq 0\},\forall y\in M.$$
By the monotonicity of the operator L and Lemma 3.1, $$P(y)\subset
J(y), \forall y\in M.$$ Since P is a $KKM$ mapping, $Co(\{x_0\ldots
x_n\})\subset\cup_{i=0}^n P(x_i)\subset \cup_{i=0}^n J(x_i)$. hence
J is $KKM$ mapping.

Furthermore, it can be easy to prove that $J(y)$ is weakly compact.
In fact, Since M is a bounded closed convex subset, it can be
convinced that M is weakly compact. But $J(y)\subset M$, thus we
only need to prove $J(y)$ is weakly closed. In order to do that let
$y_n \in J(y)\subset M$, s.t.
 $$y_n \rightharpoonup y_0 ,so ,y_n-y \rightharpoonup y_0-y,\forall y\in M.$$
According to the semi-continuity
 of L,$<L(y),y_n-y>\rightarrow <L(y),y_0-y>$. By $y_n \in J(y)$,we have $<L(y),y_n-y>\leq0$, so $<L(y),y_0-y>\leq0$. Hence $y_0\in J(y)$, that is to say, $J(y)$ is weakly closed.

From the above proof, J is a $KKM$ mapping, and $\forall y,J(y)$ is
weakly compact. According to $FKKM$ theorem, $\cap_{y\in M}
J(y)\neq\varnothing$. According to Lemma 3.1,$\cap_{y\in M} J(y) =
\cap_{y\in M} P(y)\neq\varnothing$,thus
\[\exists x_0\in \cap_{y\in M} P(y),s.t.,
\\<L(x_0),x_0-y>\leq0,\forall y\in M.\]

\section{\textbf{The Proof of Fixed Point Theorem1.3 for the Non-Expansive Mapping}}\label{sect 4}

\paragraph{}\textbf{Proof}: Let $L(x)=x-A(x),\forall x \in M $, then $L:M\rightarrow H$.By the non-expansive property for the operator A ,we will prove L= I - A is a monotonous operator:
 \[<L(x)-L(y),x-y>\geqslant0, \forall x,y \in M .\] In fact,
 \begin{eqnarray*}
<L(x)-L(y),x-y> &=& <x-y-(Ax-Ay),x-y> \\
&=&\|x-y\|^2-<Ax-Ay,x-y>\\
&\geq&\|x-y\|^2-\|Ax-Ay\|^2\cdot\|x-y\|^2\\
&\geq&\|x-y\|^2-\|x-y\|^2=0. \end{eqnarray*} Since A is
non-expansive mapping :$\forall x,y \in M$ ,$\|Ax-Ay\|\leq
\|x-y\|$,it can be derived that L is semi-continuous mapping.
According to Lemma 3.3, $\exists x_0 \in M$ ,s.t.
$$<L(x_0),x_0-x>\leq0, \forall x \in M \eqno(4.1)$$
\\\quad\quad Specially we let $x=A(x_0)\in M$ and substitute it into (4.1) ,we get:
$$<L(x_0),x_0-A(x_0)>\leq 0\eqno(4.2)$$
Since $L(x_0)=x_0-A(x_0)$, so from (4.2) we have
$\|x_0-A(x_0)\|\leq0$. It's clearly that norm $ \|x_0-A(x_0)\|\geq
0$, hence $\|x_0-A(x_0)\|=0$, that is to say, $A(x_0)=x_0.$

\section{\textbf{The Proof of  Theorem1.4 }}\label{sect 5}

\paragraph{}\textbf{Proof}: Let $(Au)(x)=\lambda\int_{a}^{b} F(x,y,u(y)) dy+f(x)$,
\begin{eqnarray*}
\|Au-Av\|^{2}&=&\int_{a}^{b} |\lambda\int_{a}^{b} K(x,y)((u(y)-v(y))
dy|^{2} dx\\
&\leq&\lambda^{2}\int_{a}^{b} [\int_{a}^{b} |K|^{2} dy \int_{a}^{b}
|u(y)-v(y)|^{2} dy] dx\\
&=&[\lambda^{2}\int_{a}^{b}\int_{a}^{b}|K(x,y)|^{2}dxdy]\cdot\|u-v\|_{L^{2}}^{2}.
\end{eqnarray*}

If  $[\lambda^{2}\int_{a}^{b}\int_{a}^{b}|K(x,y)|^{2}dxdy]\leq1$,
then A is non-expansive mapping.

Furthermore,we restrict A on a bounded closed convex subset M of
$L^{2}$ such that $AM \subseteq M$.Given $r>o$,let$ M=\{u \in L^{2}
| \|u\|_{L^{2}}\leq r\}$,
\begin{eqnarray*}
\|Au(x)\|_{L^{2}}&\leq&\|\lambda\int_{a}^{b}F(x,y,u(y))dy\|_{L^{2}}+\|f\|_{L^{2}}\\
&\leq&|\lambda|[\int_{a}^{b}|\int_{a}^{b}F
dy|^{2}dx]^{1/2}+\|f\|_{L^{2}}\\
&\leq&|\lambda|\{\int_{a}^{b}[\int_{a}^{b}|K(x,y)u(y)|dy]^{2}dx\}^{1/2}+\|f\|_{L^{2}}\\
&\leq&|\lambda|\{\int_{a}^{b}[\int_{a}^{b}|K|^{2}dy\int_{a}^{b}|u(y)|^{2}dy]dx\}^{1/2}\\
&=&|\lambda|[\int_{a}^{b}\int_{a}^{b}|K(x,y)|^{2}dxdy]^{1/2}\cdot\|u\|_{L^{2}}+\|f\|_{L^{2}}\\
&\leq&|\lambda|[\int_{a}^{b}\int_{a}^{b}|K(x,y)|^{2}dx
dy]^{1/2}\cdot r+\|f\|_{L^{2}}\\
&\leq&r.
\end{eqnarray*}

(1).If $f\equiv0$,and
$$|\lambda|[\int_{a}^{b}\int_{a}^{b}|K(x,y)|^{2}dx dy]^{1/2}\leq
1.$$
Then we can choose any given $r>0$.

(2).If $f\neq0$,and
$$|\lambda|[\int_{a}^{b}\int_{a}^{b}|K(x,y)|^{2}dx dy]^{1/2}<1.$$
Then choose r such that, $$r\geq
\frac{\|f\|_{L^{2}}}{1-|\lambda|[\int_{a}^{b}\int_{a}^{b}|K(x,y)|^{2}dx
dy]^{1/2}}.$$

{\bf{Acknowledgements.}}  I would like to thank Professor Zhang
Shiqing for his lectures on Functional Analysis ,the paper was
written when I attended his lectures, I would like to thank
Professor Shiqing Zhang for his many helpful
discussions,encouragements and suggestions and corrections,this
paper is supported partially by NSF of China.

\end{document}